\documentclass{amsart}
\usepackage{amssymb}
\newtheorem{thm}{Theorem}[section]
\newtheorem{conj}{Question}[section]
\newtheorem{lemma}[thm]{Lemma}
\newtheorem{prop}[thm]{Proposition}

\newtheorem{defn}[thm]{Definition}
\theoremstyle{definition}

\def \mrn {{\mathbb R}^n}

\def \mr {{\mathbb R}}

\def \mcs {{\mathcal S}}

\def \mcr {{\mathcal R}}

\def \mn {{\mathbb N}}

\def \eps {\epsilon}
\def \la {\lambda}
\def \ga {\gamma}

\def \lan {\langle}
\def \ran {\rangle}

\def\qed{\hfill$\square$}

\def \dive {\operatorname{div}}

\newcommand{\vol}{\operatorname{vol}}

\def \ha{ {\frac{1}{2}}}
\def \tha{ {\frac{3}{2}}}
\def \fha{ {\frac{5}{2}}}
\def \oq {\frac{1}{4}}

\def \p {\partial}

\def \rao#1 {\frac{\p}{\p #1} #1}

\numberwithin{equation}{section}
\title[A Support Theorem For The Radiation Fields]{A Support Theorem For The Radiation Fields On   Asymptotically  Euclidean Manifolds}
\author[S\'{a} Barreto]{Ant\^{o}nio S\'{a} Barreto}
\address{Department of Mathematics, Purdue University,
150 North University Street, West Lafayette IN  47907, USA}
\email{sabarre@math.purdue.edu}
\begin{document}
\input epsf
\maketitle
\section{Introduction}

We prove a support theorem for the radiation fields on  asymptotically
Euclidean manifolds with metrics which are warped products near infinity.  It generalizes to this setting the well known support theorem for the Radon transform in $\mrn.$     The main reason we are interested in proving such a theorem is the possible application to the problem of reconstructing an asymptotically Euclidean manifold from the scattering matrix at all energies, see \cite{inv}.

An asymptotically Euclidean manifold \cite{mel1,jsb} is
a $C^\infty$ compact manifold $X$  with
boundary $\p X,$ which is equipped with a $C^\infty$ Riemannian
metric $g$ that in a collar neighborhood of the boundary $\p X$
satisfies
\begin{gather}
g= \frac{dx^2}{x^4} + \frac{h(x)}{x^2},  \;\  \text{ in } \;\ [0,\eps) \times \p X, \label{metric}
\end{gather}
where $x$ is a defining function of $\p X$ and $h$ is a $C^\infty$
one parameter family of metrics on $\p X.$
The basic example is the radial compactification of $\mrn,$
\cite{mel1}.

 In this paper will  consider the class of metrics $g$ which have the following special form near
 $\p X:$
\begin{gather}
g=\frac{dx^2}{x^4}+ \psi(x) \frac{h_0}{x^2}, \;\ x \in [0,\eps), \label{warped}
\end{gather}
where $\psi\in C^\infty([0,\eps)),$ $\psi(x)>0,$ $\psi(0)=1,$ and $h_0$ is a $C^\infty$ metric on $\p X.$   
These are known as warped product metrics \cite{megs}.

We consider the wave equation on $X.$ Let  $\Delta_g$ be the Laplace operator on $X,$ and let
 $u(t,z)$ satisfy
\begin{gather}
\begin{gathered}
(D_t^2-\Delta_g)u(t,z)=0 \text{ on } \mr \times X, \\
u(0)=f_1, \;\ D_t u(0)=f_2, \;\ \text{ with } \;\  f_1, f_2 \in C_0^\infty(X).
\end{gathered}\label{weq}
\end{gather}
A function $f \in C_0^\infty(X)$ if it is $C^\infty$ and its support does not intersect the boundary of $X.$

The following is proved in \cite{fried1,fried2}:
\begin{thm}\label{thm1}  Let $x$ be the boundary defining function for which \eqref{metric} holds, and let 
$z=(x,y),$ $y\in \p X,$ be the corresponding boundary normal coordinates in a collar neighborhood of the boundary.  Then
\begin{gather}
\begin{gathered}
v_+(x,s,y)=x^{-\frac{n-1}{2}} u( s+ \frac{1}{x}, x,y) \in C^\infty(\mr_s\times [0,\eps)_x \times \p X), \\
v_-(x,s,y)=x^{-\frac{n-1}{2}} u( s- \frac{1}{x}, x,y) \in C^\infty(\mr_{s}\times [0,\eps)_x \times \p X).
\end{gathered} \label{radf}
\end{gather}
\end{thm}
Friedlander \cite{fried1,fried2} defined the forward and backward radiation fields respectively as
\begin{gather}
\begin{gathered}
\mcr_+(f_1,f_2)= D_s v_+(0,s,y) \;\ \text{ and } \;\
\mcr_-(f_1,f_2)= D_s v_-(0,s,y).
\end{gathered} \label{mcr+}
\end{gather}
Lax and Phillips  \cite{lp} proved that in $\mrn$  the forward (or backward) radiation field  is
the modified Radon transform, that is:
\begin{gather*}
 \mcr_{+}(f_1,f_2)(s,\omega)= |D_s|^{\frac{n-3}{2}} R f_1(s,\omega)+
 |D_s|^{\frac{n-1}{2}} R f_2(s,\omega), \text{ where } \\
 R f(s,\omega)=\int_{\lan x,\omega\ran=s} f(z) \; d\sigma, \;\ \sigma \text{ 
 is the surface measure on } \lan x,\omega\ran=s,  \text{ is the Radon transform.} 
\end{gather*}

Helgason's celebrated support theorem for Radon transforms \cite{helg} says that if $f$ is a rapidly decaying function in $\mrn$ and $R f(s,\omega)=0$ for $s<s_0,$ $s_0<0$ (and hence by symmetry $Rf(s,\omega)=0$ for $|s|>|s_0|$) then $f$ is supported in the ball of radius $|s_0|.$  
The assumption that $f$ is rapidly decaying cannot be entirely removed.
For example,  for any $m\in \mn,$ there are smooth functions  $f(z)$ in $z\in \mrn$
which are not compactly supported, decay like $|z|^{-m},$ and whose Radon transform is compactly supported. See for example \cite{helg}.  

This is a result in control theory, where the support of a function can be exactly controlled by the support of its Radon transform.  We want to address the analogue question for Radiation fields. 
The following was proved in \cite{euradf}:
\begin{thm} \label{compsup}If $f\in C_0^\infty(X),$ $g$ is an arbitrary asymptotically Euclidean metric, and $\mcr(0,f)(s,y)=0$ for $s<-\frac{1}{x_0},$ $x_0\in (0,\eps),$ then
$f=0$ if $x<x_0.$
\end{thm}

This says that if there exists some $x_1\in (0,\eps),$ such that $f(x,y)=0$ for $x<x_1,$ but $\mcr_+(0,f)(s,y)=0$ for $s<-\frac{1}{x_0},$ and $x_0>x_1,$
then in fact $f(x,y)=0$ if $x<x_0.$  The purpose of this paper is to discuss the following
\begin{conj}\label{conjecture} Let $(X,g)$ be an asymptotically Euclidean manifold, and let $\mcs(X)$ be the space of functions in $C^\infty(X)$
which are smooth up to $\p X$ and vanish to infinite order at $\p X.$
If $f\in \mcs(X)$ and $\mcr_+(0,f)(s,y)=0$ for $s<-\frac{1}{x_0},$ $x_0\in (0,\eps),$ is it true that
$f=0$ if $x<x_0?$ 
\end{conj}

 We answer this question in the affirmative in the following particular case:
\begin{thm}\label{main}   If $g$ is a warped product metric, the dimension of $X$ is greater than or equal to 3,  $f\in \mcs(X),$ and $\mcr_+(0,f)(s,y)=0$ for $s<-\frac{1}{x_0},$ then $f=0$ if $x<x_0.$
\end{thm}

\section{Energy Estimates}

The first step in the proof of Theorem \ref{main} is to obtain 
estimates  the solution to \eqref{weq} up to $x=0$ and $s=-\infty.$

The Laplacian with respect to the metric \eqref{warped} is, in a
neighborhood of $\p X,$ given by
\begin{gather*}
\Delta_g= -x^4 \p_x^2 +(n-3)x^3 \p_x - x^4 A(x)\p_x + x^2\psi(x)^{-1} \Delta_{h_0},
\end{gather*}
where $A(x)=\p_x \log \left(\psi(x)^\frac{n-1}{2}\right),$  and
$\Delta_{h_0}$ is the Laplacian on $\p X$ with respect to the metric
$h_0.$ In what follows it is convenient to get rid of first order
terms, so we will work with $Q=F(x)^{-1}\Delta_g F(x),$ with
$F(x)=x^{\frac{n-1}{2}}\left(\psi(x)\right)^{-\frac{n-1}{4}}.$ We
get that
\begin{gather*}
Q=-(x^2\p_x)^2 +x^2 \phi(x)\Delta_{h_0}+x^2 B(x), \text{ where } \\
\phi(x)=[\psi(x)]^{-1}, \;\ B(x)=\frac{(n-1)(n-3)}{4} +x B_1(x), \;\ B_1\in
C^\infty([0,\eps)).
\end{gather*}
The wave equation \eqref{weq} is translated into
\begin{gather}
\begin{gathered}
(\p_t^2-(x^2\p_x)^2 +x^2 \phi(x)\Delta_{h_0}+x^2 B(x))u=0, \\
u(0)=f_1(x,y), \;\ \p_tu(0)=f_2(x,y).
\end{gathered}\label{weqnew}
\end{gather}

Instead of working with coordinates $x$ and $s,$
and the forward and backward radiation fields separately,  it is better to  work with the forward and backward
radiation fields simultaneously.
So we define
\begin{gather}
s_+= t-\frac{1}{x} \text{ and } s_-= t+\frac{1}{x}, \;\ x>0. \label{rw}
\end{gather}

Let $u$ be the solution to the wave equation
\eqref{weqnew}, and let $w=F(x)u.$  In coordinates \eqref{rw},
$$w(s_+,s_-,y)=F(\frac{2}{s_--s_+}) u(\frac{s_++s_-}{2},
\frac{2}{s_--s_+}, y)$$ satisfies
\begin{gather}
\begin{gathered}
(s_- -s_+)^2 \p_{s_+} \p_{s_-}w+
\phi\left(\frac{2}{s_--s_+}\right) \Delta_{h_0}w+ B\left(\frac{2}{s_--s_+}\right) w=0, \\
w(s_+,-s_+,y)= F(-\frac{1}{s_+}, y) f_1(-\frac{1}{s_+}, y), \\ 
(\p_{s_+} w)(s_+,-s_+,y)= \ha F(-\frac{1}{s_+})f_2(-\frac{1}{s_+}, y) +  
\frac{1}{2s_+^2} F'(-\frac{1}{s_+})f_1(-\frac{1}{s_+}, y)+
\frac{1}{2s_+^2} F(-\frac{1}{s_+})f_1'(-\frac{1}{s_+}, y),
\end{gathered}\label{weq2}
\end{gather}

 We want to understand the behavior of $w$ as $s_+\sim
-\infty,$ and $s_-\sim\infty$ and thus we compactify $\mr \times X$ by
setting
\begin{gather*}
\mu= -\frac{1}{s_+}  \text{ and } \nu= \frac{1}{s_-}
\end{gather*}

So 
\begin{gather}
w(\mu,\nu,y)=F\left(\frac{2\mu\nu}{\mu+\nu}\right) u\left( \frac{\mu-\nu}{\mu+\nu}, \frac{2\mu\nu}{\mu+\nu},y\right) \label{wfu}
\end{gather}
satisfies
\begin{gather}
\begin{gathered}
\left((\mu+\nu)^2 \p_{\mu} \p_{\nu} -\phi\left(\frac{2\mu\nu}{\mu+\nu}\right)\Delta_{h_0}w
- B \left(\frac{2\mu\nu}{\mu+\nu}\right)\right)w=0, \\
w(\mu,\mu,y)=\widetilde{f_1}(\mu,y), \;\ (\p_\mu
w)(\mu,\mu,y)=\widetilde{f_2}(\mu, y),
\end{gathered}\label{weq3}
\end{gather}
where 
\begin{gather*}
\widetilde{f_1}(\mu,y)= F(\mu) f_1(\mu,y), \text{ and } \\
\widetilde{f_2}(\mu,y)=\frac{1}{2\mu^2} F(\mu)f_2(\mu,y)+ \ha F'(\mu)f_1(\mu,y)+ \ha F(\mu)f_1'(\mu,y).
\end{gather*}

Since we are dealing with a degenerate equation, we will work with weighted Sobolev spaces.   
\begin{defn}\label{definition}  Let $T>0$ and $\Omega_T=(0,T)\times (0,T).$  We define
\begin{gather*}
H^s_{j}([0,T]\times \p X)=\{ f\in L^2([0,T]\times \p X), \;\  \mu^{-j} f \in H^s([0,T]\times \p X)\}, \text{ with norm } ||f||_{s,j}, \text{ and } \\
H^s_j(\Omega_T \times \p X)=\{f\in L^2(\Omega_T \times \p X), \;\  
(\mu+\nu)^{-j} f \in H^s(\Omega_T\times \p X)\}, 
\text{ with norm } |||f|||_{s,j}, \text{ and } \\
\end{gather*}
\end{defn}

The next step is to prove
\begin{thm}\label{smooth} Let  $\widetilde{f_j}(\mu,y)\in C^\infty([0,T] \times \p X),$ $j=1,2,$ be such that $\p_{\mu}^k \widetilde{f_j}(0,y)=0,$
 $k=0,1,2,....$ Let  $w$ satisfy \eqref{weq3} in $\Omega_T\times \p X.$  Then there exists
 $T_0>0$ such that $w$ has a $C^\infty$ extension up to
$\overline{\Omega_T}  \times \p X,$ $T\leq T_0.$
\end{thm}
\begin{proof}   

By finite speed of propagation, $w\in C^\infty(\Omega_T\times \p X),$ we want to establish the regularity up to the closure $\overline{\Omega_T}\times \p X.$
We begin the proof with the following elementary lemmas:

\begin{lemma}\label{bound0}  Let $\mu \in (0,T),$ and $w \in C^\infty(\Omega_T).$ Then the following inequalities are true
\begin{gather}
\begin{gathered}
\int_\mu^b (\mu+\nu)^{-1-k}|w(\mu,\nu)|^2 \; d\nu \leq 2b \mu^{-1-k} |w(\mu,\mu)|^2 +b \int_\mu^b
(\mu+\nu)^{-k} |\p_\nu w(\mu,\nu)|^2 \; d\nu, \;\ k\in \mn,
\end{gathered}\label{id1gen}
\end{gather}
and
\begin{gather}
\begin{gathered}
\int_\mu^\nu (\mu+\nu)^{-1}|w(\mu,\nu)|^2 \; d\mu \leq 2 \nu^{-1}|w(\nu,\nu)|^2 +\nu \int_a^\nu
 |\p_\mu w(\mu,\nu)|^2 \; d\mu, \;\ k\in \mn.
\end{gathered}\label{id2gen}
\end{gather}
\end{lemma}
\begin{proof} We  write
\begin{gather}
w(\mu,\nu)=w(\mu,\mu)+\int_\mu^\nu \p_s w(\mu,s) \; ds. \label{ftc}
\end{gather}
The Cauchy-Schwartz inequality gives
\begin{gather*}
|w(\mu,\nu)|^2\leq 2 |w(\mu,\mu)|^2+2(\nu-\mu) \int_\mu^\nu |\p_s
w(\mu,s)|^2 \; ds.
\end{gather*}
Since $\nu>0$ 
\begin{gather*}
\int_\mu^b (\mu+\nu)^{-1-k}|w(\mu,\nu)|^2 \; d\nu \leq 2(b-\mu) \mu^{-1-k} |w(\mu,\mu)|^2+2\int_\mu^b
\int_\mu^\nu (\mu+\nu)^{-k} |\p_s w(\mu,s)|^2 \; ds d\nu\leq \\
2b \mu^{-1-k}|w(\mu,\mu)|^2  + b\int_\mu^b (\mu+\nu)^{-k}|\p_s w(\mu,s)|^2 \; d\nu.
\end{gather*}
The proof of \eqref{id2gen} is identical, and this ends the proof of the Lemma.
\end{proof}

\begin{lemma}\label{intFlemma} Let $\Omega_{a,b}\subset \Omega_T$ be the region
defined by
\begin{gather}
\Omega_{a,b}=\{(\mu,\nu): \nu\geq \mu\geq 0, \;\ 0 \leq a \leq \mu, \;\ \nu\leq b\}, \label{oab}
\end{gather}
 let $F(\mu,\nu)\in L^1(\Omega_T)$ and let $a_0\leq a \leq b.$ Then

\begin{gather}
\int_{a_0}^b \left(\int_{\Omega_{ab}} F(\mu,\nu) \; d\mu d\nu \right) \; da =
\int_{\Omega_{a_0 b}} (\mu-a_0) F(\mu,\nu) \; d\mu d\nu. \label{intF}
\end{gather}
\end{lemma}
The proof is a straightforward application of Fubini's theorem.

\begin{lemma}\label{decaym}  If    $\int_0^{T} \mu^{-m-2}|w(\mu,\mu)|^2 d\mu$ and
$\int_{\Omega_{0T}} (\mu+\nu)^{-m} |(\p_\nu-\p_\mu)w|^2  d\mu d\nu$ are finite, then
\begin{gather}
\begin{gathered}
\int_0^T\int_0^T (\mu+\nu)^{-m-2} |w|^2  d\mu d\nu \leq  \\
2^{-m}  \int_0^{T}  \mu^{-m-1} w(\mu,\mu) \; d\mu +
 \int_0^T\int_0^T   (\mu+\nu)^{-m} |(\p_\nu-\p_\mu)w|^2  d\mu d\nu. 
\end{gathered}\label{decay5}
\end{gather}
\end{lemma}
\begin{proof}
To see this it is better to rotate the axes $\mu$ and $\nu$ and use coordinates
\begin{gather}
\begin{gathered}
r=\mu+\nu , \;\  \tau=\nu-\mu, \text{ so } \\
\mu=\ha(r-\tau), \;\ \nu=\ha(r+\tau).
\end{gathered}  \label{rotate}
\end{gather}
In the region $\mu\geq 0,$ and $T\geq \nu\geq \mu,$ we have
\begin{gather*}
T\sqrt{2} \geq r\geq \tau, \;\  \tau\geq 0.
\end{gather*}
In these coordinates $\p_\nu-\p_\mu=2\p_\tau$ and the diagonal $\mu=\nu$ becomes $\tau=0.$
We write, for $r\geq \tau,$
\begin{gather*}
w (\tau, r)= w(0,r) + \int_0^\tau \p_s w(s,r) \; ds
\end{gather*}
Using the Cauchy-Schwartz inequality and the fact that $\tau/r\leq 1,$ we have
\begin{gather*}
r^{-m-2} |w(\tau,r)|^2  \leq 2 r^{-m-2} |w_k(0,r)|^2 + 2  \tau r^{-m-2}
\int_0^\tau  | \p_s w_k(s,r)|^2 \; ds \leq \\
2 r^{-m-2} |w_k(0,r)|^2 + 2  r^{-m-1}
\int_0^\tau  | \p_s w_k(s,r)|^2 \; ds 
\end{gather*}
Therefore
\begin{gather*}
\int_0^{T\sqrt{2}} \int_0^r r^{-m-2} |w^2(\tau,r)|^2 \; d\tau dr \leq 
2 \int_0^{T\sqrt{2}} r^{-m-1} |w(0,r)|^2 \; dr + 2 \int_0^{T\sqrt{2}} \int_0^r r^{-m-1} \int_0^\tau  | \p_s w_k(s,r)|^2 \; ds\; d\tau dr
= \\ 2 \int_0^{T\sqrt{2}} r^{-m-2} |w(0,r)|^2 \; dr + 
2 \int_0^{T\sqrt{2}} \int_0^r r^{-m-1} (r-s)  | \p_s w(s,r)|^2 \; ds\;  dr.
\end{gather*}
We get the same bound in the region $r \geq -\tau.$ Translating this back into coordinates $\mu$ and $\nu$ we get \eqref{decay5}.
\end{proof}

Now we prove uniform energy estimates up to $\{\mu=0\},$ $\{\nu=0\}.$  When $n>3$ one can choose
$T_0$ such that if $\mu<T_0$ and $\nu<T_0,$ $B> \frac{(n-1)(n-3)}{8}.$ When $n=3,$ $B$ is not necessarily positive.   In this case it is convenient to work with an eigenfunction decomposition.
Let $\phi_k,$ $k\in \mn,$ be the eigenfunctions of $\Delta_{h_0}$ and let $\{\lambda_k\},$ $k\in \mn,$
with $0=\lambda_{k-1}\leq \lambda_k$ be the corresponding eigenvalues.  Let $w$ be the solution to \eqref{weq3}  and let
$w_k=\lan w, \phi_k\ran_{L^2(\p X)}.$ When $n=3$ and $k=1,$ $w_1$ satisfies
\begin{gather}
\begin{gathered}
\left((\mu+\nu)^2\p_\mu\p_\nu-  \frac{2\mu\nu}{\mu+\nu} B_1(\frac{2\mu\nu}{\mu+\nu})\right)w_1=0 
\text{ in } (0,T)\times(0,T), \\
w(\mu,\mu)=q_1(\mu), \;\ \p_\mu w(\mu,\mu)=q_2(\mu), \;\ 
\end{gathered}\label{degen}
\end{gather}
with $q_1=\lan \widetilde{f_1}, \phi_1\ran$ and  $q_2=\lan \widetilde{f_2}, \phi_1\ran.$ It is clear that
$||q_1||_{s,j}<\infty$ and $||q_2||_{s,j}<\infty$ for every $s$ and $j.$

 We have
\begin{prop}\label{deg1}  Let  $\Omega_T=(0,T)\times (0,T)$ and let $w_1\in C^\infty(\Omega_T)$ satisfy \eqref{degen} in
$\Omega_T.$  Suppose that $||q_1||_{s,j}<\infty$ and $||q_2||_{s,j}<\infty,$ for all $s$ and $j.$
 Then there exists $T_0>0$ such that $w_1\in C^\infty(\overline{\Omega_T}),$ $T\leq T_0.$
\end{prop}
\begin{proof} The proof relies on the following energy estimates:
\begin{lemma} Let  $w_1\in C^\infty(\Omega_T)$ satisfy \eqref{degen} in $\Omega_T.$  There
exist $T_0$ and a constant $C=C(T_0)>0$ such that if $T\leq T_0,$
\begin{gather}
\begin{gathered} 
\text{ for any fixed  } \mu \in [0,T], \;\
\int_0^T  |\p_\nu w_1(\mu,\nu)|^2  d\nu \leq C(||q_1||_{1,\ha}^2+ ||q_2||_{1,\ha}^2), \\
\text{  and  for any fixed  } \nu \in [0,T], \;\ 
\int_0^T  |\p_\nu w_1(\mu,\nu)|^2  d\mu \leq C(||q_1||_{1,\ha}^2+ ||q_2||_{1,\ha}^2).\\
\end{gathered} \label{mainest2deg} 
\end{gather}
\end{lemma}
\begin{proof}
To prove this we first multiply \eqref{degen} by $(\mu+\nu)^{-2}(\p_\mu-\p_\nu)w_1.$
 We obtain
\begin{gather*}
\ha \p_{\nu}\left( (\p_\mu w_1)^2 + 2\mu\nu (\mu+\nu)^{-3} B_1 w_1^2 \right)
-\ha \p_{\nu}\left( (\p_\nu w_1)^2 +2 \mu\nu (\mu+\nu)^{-3} B_1 w_1^2 \right)
-\ha(\mu+\nu)^{-2} B_2 w_1^2=0,
\end{gather*}
where $B_2= 2\frac{\mu-\nu}{\mu+\nu} B_1+ 4\frac{\mu\nu(\mu-\nu)}{(\mu+\nu)^2} B_1'(\frac{2\mu\nu}{\mu+\nu}).$ After we integrate it in $\Omega_{ab}$ we obtain,
\begin{gather*}
\ha \int_{a}^b  \left( (\p_\nu w_1)^2 + 2\mu\nu (\mu+\nu)^{-3} B_1 |w_1|^2 \right)(a,\nu) d\nu +
\ha \int_{a}^b  \left( (\p_\mu w_1)^2 + 2\mu\nu (\mu+\nu)^{-3} B_1 |w_1|^2 \right)(\mu,b) d\mu 
- \\ \ha \int_{\Omega_{ab}} (\mu+\nu)^{-2} B_2 |w_1|^2 d\mu d\nu= 
\frac{1}{2\sqrt{2}}\int_a^b \left(( q_1'+q_2)^2 + q_2^2 + \ha \mu^{-1}B_1(\mu) q_1^2\right) d\mu.
\end{gather*}

We apply Lemma \ref{bound0}  to show that we can pick $T_0$ such that for $T<T_0,$ 
\begin{gather*}
 \int_a^b 2 \mu\nu (\mu+\nu)^{-3} |B_1| |w_1|^2(a,\nu)  d\nu \leq
 2 a^{-1} q_1^2(a) + \oq\int_{a}^b  (\p_\nu w_1)^2 (a,\nu) d\nu, \text{ and }\\
  \int_{a}^b 2 \mu\nu (\mu+\nu)^{-3} |B_1| |w_1|^2(\mu,b)  d\mu \leq
 2 b^{-1} q_1^2(b) + \oq \int_{a}^b  (\p_\nu w_1)^2(\mu,b) d\mu.
 \end{gather*}
 Therefore
 \begin{gather}
  \begin{gathered}
 \oq \int_a^b   (\p_\nu w_1)^2(a,\mu) d\nu + \oq \int_a^b  (\p_\mu w_1)^2(\mu,b)  d\mu 
- \ha \int_{\Omega_{ab}} (\mu+\nu)^{-2} |B_2 ||w_1|^2 d\mu d\nu\leq \\
a^{-1} q_1^2(a)+b^{-1} q_1^2(b)+  \frac{1}{2\sqrt{2}}\int_a^b\left( (q_1'+q_2)^2 + q_2^2 + \ha \mu^{-1}B_1(\mu) q_1^2\right) d\mu.
\end{gathered}\label{imp}
\end{gather}
If we drop the second integral from this inequality and integrate the remaining terms in $a,$ with
$a_0\leq a \leq b,$ and use \eqref{intF} we get that
\begin{gather*}
 \oq \int_{\Omega_{a_0b}}   (\p_\nu w_1)^2 d\nu -
  \int_{\Omega_{a_0b}} \mu(\mu+\nu)^{-2} |B_2| |w_1|^2 d\mu d\nu\leq \\
 \int_{a_0}^b \mu^{-1} q_1^2(\mu) d\mu +   q_1^2(b) +
\frac{ T}{2\sqrt{2}} \int_{a_0}^b \left( (q_1'+q_2)^2 + q_2^2 +  q_1^2\right) d\mu.
\end{gather*}
We can use Lemma \ref{bound0} to show that 
\begin{gather}
 \int_{\Omega_{a_0b}} (\mu+\nu)^{-1} |w_1|^2 d\mu d\nu\leq
T \int_{\Omega_{a_0b}}   (\p_\nu w_1)^2 d\nu  + \int_{a_0}^b \mu^{-1} q_1^2(\mu) \; d\mu. \label{imp1}
 \end{gather}
and therefore, if $T_0$ is small,
 \begin{gather}
 \frac{1}{8} \int_{\Omega_{a_0b}}   (\p_\nu w_1)^2 d\nu \leq q_1^2(b) +
  C(||q_1||_{1,\ha}^2+ ||q_2||_{1,\ha}^2). \label{imp2} 
\end{gather}
Now we substitute \eqref{imp2} and \eqref{imp1} into \eqref{imp} and use that
$|| q_1||_{L^\infty} \leq || q_1||_{1,0}$ and deduce that
\begin{gather*}
 \oq \int_{a}^b   (\p_\nu w_1)^2(a,\nu) d\nu + \oq \int_a^b (\p_\mu w_1)^2(\mu,b)  d\mu 
\leq C (||q_1||_{1,\ha}^2+ ||q_2||_{1,\ha}^2). 
\end{gather*}
 By symmetry this estimate also holds in the region below the diagonal.  This proves
  \eqref{mainest2deg}.  
\end{proof}
Now we prove Proposition \ref{deg1}.   For $\nu\geq \mu$ we write
\begin{gather*}
w_1(\mu,\nu) =w_1(\mu,\mu) + \int_\mu^\nu \p_s w_1(\mu,s) ds.
\end{gather*}
By the Cauchy-Schwartz inequality,
\begin{gather*}
|w_1(\mu,\nu)|^2 \leq 2|w_1(\mu,\mu)|^2 + 2(\nu-\mu) \int_\mu^\nu |\p_s w_1(\mu,s)|^2 ds,
\end{gather*}
and we deduce from  \eqref{mainest2deg} that if $(\mu,\nu)\in \Omega_T,$ with $T<T_0$ small,
\begin{gather*}
(\mu+\nu)^{-1}|w_1(\mu,\nu)|^2 \leq 2 \mu^{-1}q_1^2(\mu) + 2 \int_\mu^\nu |\p_s w_1(\mu,s)|^2 ds
\leq C( ||q_1||_{1,\ha}^2+ ||q_2||_{1,\ha}^2).
\end{gather*}
By symmetry with respect to the diagonal,
\begin{gather*}
|w_1(\mu,\nu)|\leq C(  (||q_1||_{1,\ha}^2+ ||q_2||_{1,\ha}^2)^\ha (\mu+\nu)^{\ha}, \text{ in } \Omega_{T_0}.
\end{gather*}
From now on we will use $C(q_1,q_2)$ to denote a constant which depends on
the norms $||q_1||_{s,j}$ and $||q_2||_{s,j}$ for some $s$ and $j.$ 
We then go back to equation \eqref{degen} and deduce that if $(\mu,\nu)\in \Omega_{T_0},$ 
\begin{gather*}
|\p_\mu\p_\nu w_1|= |2\mu\nu(\mu+\nu)^{-3} B_1 w_1| \leq C(q_1,q_2) (\mu+\nu)^{-\ha}
\end{gather*}
thus
\begin{gather*}
|\p_\nu w_1(\mu,\nu)|\leq |\p_\nu w_1(\nu,\nu)|+ \int_\mu^\nu |\p_s \p_\nu w_1(s,\nu) | ds
\leq  \\  |\p_\nu w_1(\nu,\nu)|+ C(q_1,q_2)\int_\mu^\nu (\mu+s)^{-\ha}  ds \leq
C(q_1,q_2) (\mu+\nu)^{\ha} .
\end{gather*}
A similar  argument shows that
\begin{gather*}
|\p_\mu w_1(\mu,\nu)| \leq C(q_1,q_q) (\mu+\nu)^{\ha}, \text{ if } \nu\geq \mu.
\end{gather*}
By symmetry these estimates hold below the diagonal. This implies that
\begin{gather*}
|w_1(\mu,\nu)| \leq C(q_1,q_2) (\mu+\nu)^{\tha}  \text{ in } \Omega_{T_0}.
\end{gather*}
We then differentiate  equation \eqref{degen} and find that
\begin{gather*}
\p_\mu \p_\nu^2 w_1=\left[\mu(2\mu-\nu)(\mu+\nu)^{-4}B_1 -2\mu^3\nu(\mu+\nu)^{-5} B_1' \right] w_1
+ 2\mu\nu(\mu+\nu)^{-1} B_1 \p_\nu w_1.
\end{gather*}
We deduce that
\begin{gather*}
|\p_\mu \p_\nu^2 w_1| \leq C(q_1,q_2)(\mu+\nu)^{-\ha}
\end{gather*}
which implies that
\begin{gather*}
|\p_\mu \p_\nu w_1(\mu,\nu)| \leq C(q_1,q_q) (\mu+\nu)^{\ha}  \;\
|\p_\nu^2 w_1(\mu,\nu)| \leq C(q_1,q_q) (\mu+\nu)^{\ha},  \;\ \text{ in } \Omega_{T_0}.
\end{gather*}
A similar argument gives that
\begin{gather*}
|\p_\mu^2 w_1(\mu,\nu)| \leq C(q_1,q_q) (\mu+\nu)^{\ha},  \;\ \text{ in } \Omega_{T_0}.
\end{gather*}
hence we deduce that
\begin{gather*}
|\p_\mu w_1(\mu,\nu)| \leq C(q_1,q_q) (\mu+\nu)^{\tha},  \;\
|\p_\nu w_1(\mu,\nu)| \leq C(q_1,q_q) (\mu+\nu)^{\tha}.
\end{gather*}

This argument can be repeated to show that 
\begin{gather*}
|\p_\mu^j \p_\nu^k w_1(\mu,\nu)|\leq C(q_1,q_2) (\mu+\nu)^m, \;\ j, k, m\in \mn  \text{ in } \Omega_{T_0}.
\end{gather*}
This implies the claim of  Proposition \ref{deg1}.
\end{proof}

Next we study the non-degenerate cases, i.e either $n>3$ or if $n=3,$  $\int_{\p X} w=0,$ i.e
$w$ is orthogonal to the first eigenfunction.

\begin{lemma}\label{trueestimates} Let  $\Omega_T=(0,T)\times (0,T)$ and let 
$W \in C^\infty(\Omega_T)$ satisfy
\begin{gather}
\begin{gathered}
\left((\mu+\nu)^2\p_\mu\p_\nu- \phi(\frac{2\mu\nu}{\mu+\nu}) \Delta_{h_0} - B(\frac{2\mu\nu}{\mu+\nu})\right)W=G(\mu,\nu)  \text{ in } \Omega_T\times \p X \\
W(\mu,\mu)=q_1(\mu), \;\ \p_\mu W(\mu,\mu)=q_2(\mu). \;\ 
\end{gathered}\label{nonzero}
\end{gather}
If $n=3$ we assume that $\int_{\p X} W =0.$ Then
there exists $T_0>0,$ depending on $B$ and $\phi,$ and a constant $C$  depending on
$B$ and $\phi$ such
that  if  $|||G|||_{0,2}<\infty,$ and $||q_j||_{1,1}<\infty,$ $j=1,2,$
\begin{gather}
\begin{gathered} 
\text{ for any fixed  } \mu \in [0,T] \\
\int_0^T \int_{\p X} ( |\p_\nu W(\mu,\nu))|^2 + (\mu+\nu)^{-2}( |\nabla_{h_0}W(\mu,\nu)|^2 +  |W(\mu,\nu|^2 ) d\nu \leq
C(||q_1||_{1,1}^2+ ||q_2||_{1,1}^2+ |||G|||_{0,2}^2), \\
\text{  and  for any fixed  } \nu \in [0,T] \\
\int_0^T \int_{\p X} ( |\p_\nu W(\mu,\nu|^2 + (\mu+\nu)^{-2}( |\nabla_{h_0}W(\mu,\nu)|^2 +  |W(\mu,\nu)|^2 ) d\mu \leq
C(||q_1||_{1,1}^2+ ||q_2||_{1,1}^2+ |||G|||_{0,2}^2).\\
\end{gathered} \label{mainest2} 
\end{gather}
If $||q_1||_{1,\tha}<\infty$ and $||q_2||_{1,\tha}<\infty,$ we also have
\begin{gather} 
\begin{gathered} 
\int_{\Omega_T \times \p X}
(\mu+\nu)^{-2}(\mu |\p_\mu W|^2 + \nu |\p_\nu W|^2) 
+  (\mu+\nu)^{-3}(  |\nabla_{h_0}W|^2 +  |W|^2 ) d\mu d\nu d\vol_{h_0}   \leq \\
  CT ( ||q_1||_{2,\tha}^2+ ||q_2||_{2,\tha}^2+ |||G|||_{0,2}^2).
\end{gathered} \label{mainest3} 
\end{gather}
If $||q_j||_{2,2} <\infty, $ $j=1,2$ and $|||G|||_{0,\frac{5}{2}}<\infty,$ then
\begin{gather}
\int_{\Omega_T \times \p X} (\mu+\nu)^{-4} (  |\nabla_{h_0}W|^2 +  |W|^2 ) d\mu d\nu d\vol_{h_0}   \leq
C|||G|||_{0,\fha}+  CT ( ||q_1||_{2,2}^2+ ||q_2||_{2,2}^2).\label{mainest} 
\end{gather}
\end{lemma}
\begin{proof}  To prove these estimates we multiply \eqref{nonzero} by $(\mu+\nu)^{-m}(\p_\mu -\p_\nu)W,$ with $m\in \mr_+$, and integrate the product in $\Omega_{ab},$ with $a\leq b \leq T.$  The product is equal  to
\begin{gather*}
\ha \p_\nu\left[ (\mu+\nu)^{2-m}|\p_\mu W|^2 + (\mu+\nu)^{-m}(\phi |\nabla_{h_0} W|^2 + B|W|^2)\right]- \\
\ha \p_\mu\left[ (\mu+\nu)^{2-m}|\p_\nu W|^2 + (\mu+\nu)^{-m}(\phi |\nabla_{h_0} W|^2 + B|W|^2)\right]+ \\
\frac{m-2}{2}(\mu+\nu)^{1-m}((\p_\mu W)^2 -(\p_\nu W)^2) +
\dive_{h_0}\left[(\mu+\nu)^{-m} \phi \nabla_{h_0} W (\p_\mu-\p_\nu) W\right] + \\
(\nu-\mu)(\mu+\nu)^{-m-1}(\phi' |\nabla_{h_0} W|^2 + B' | W|^2), 
\end{gather*}
where $\phi'$ and $B'$ denote the derivative of $\phi$ and $B.$
Integrating this in $\Omega_{ab} \times \p X$ and using the divergence theorem
and  the fact that  $(\p_\nu W)(\mu,\mu)=q_1'(\mu)-q_2(\mu),$
\begin{gather}
\begin{gathered}
\ha \int_{\Sigma_1\times \p X}\left[ (\mu+\nu)^{2-m} |\p_\nu W|^2 +(\mu+\nu)^{-m}(\phi |\nabla_{h_0} W|^2 + B|W|^2)\right] d\nu d\vol_{h_0} - \\
\ha \int_{\Sigma_2 \times \p X}\left[ (\mu+\nu)^{2-m} |\p_\mu W|^2 +
(\mu+\nu)^{-m}(\phi |\nabla_{h_0} W|^2 + B|W|^2)\right] d\mu  d\vol_{h_0}  + \\
\frac{m-2}{2}\int_{\Omega_{ab} \times \p X} (\mu+\nu)^{1-m}((\p_\mu W)^2 -(\p_\nu W)^2) d\mu d\nu 
d\vol_{h_0} + \\
\int_{\Omega_{ab} \times \p X} (\nu-\mu)(\mu+\nu)^{-1-m}( \phi'|\nabla_{h_0} W|^2 + B'|W|^2) d\mu d\nu d\vol_{h_0}= \\
\int_{\Omega_{ab} \times \p X} G(\mu+\nu)^{-m} (\p_\mu-\p_\nu) W d\mu d\nu d\vol_{h_0} + \\
\frac{1}{2\sqrt{2}} \int_{\Sigma_3 \times \p X} (2\mu)^{2-m} ((q_2)^2  +
(q_2'-q_1)^2) + (2\mu)^{-m} (\phi(\mu)|\nabla_{h_0} q_1
|^2+ B(\mu) |q_1|^2) \; d\mu
d\vol_{h_0}.
\end{gathered}\label{divth}
\end{gather}
Here
\begin{gather}
\begin{gathered}
\Sigma_1=\Sigma_1(a,b)=\{(a,\nu), \; a\leq \nu \leq b\}, \;\
\Sigma_2=\Sigma_2 (a,b)=\{(\mu,b), \; 0\leq \mu \leq b\} \text{ and } \\
\Sigma_3=\Sigma_3(a,b)=\{(\mu,\nu), \; \mu=\nu, \;\ 0\leq \mu \leq b\}. 
\end{gathered}\label{defsig}
\end{gather}

When $n>3,$ and $T$ is small, $B$ is positive, but as we saw before,
when $n=3$ this is not necessarily the case.  So when $n>3$ we guarantee that the first two integrals in \eqref{divth} are positive.  When $n=3,$ $B(0)=0.$ But we assumed that in this case
$\int_{\p X} W d\vol_{h_0}=0$ and therefore
$$\int_{\p X} |\nabla_{h_0} W|^2  d\vol_{h_0}\geq \la_2 \int_{\p X} |W|^2  d\vol_{h_0}.$$
Since $\la_2>0,$ if $T_0$ is small the term in $B|W|^2$ can be absorbed by the term
in $|\nabla_{h_0} W|^2.$

So we may assume that the second integral is positive.  We drop it from \eqref{divth} and integrate the remaining terms in the variable  $a,$ which determines $\Sigma_1,$ with $a_0 \leq a \leq b.$  Using \eqref{intF} we obtain
\begin{gather}
\begin{gathered}
\ha \int_{\Omega_{a_0b}\times \p X}\left( \left[ (\mu+\nu)^{2-m}-(m-2)\mu(\mu+\nu)^{1-m}\right] |\p_\nu W|^2 + (m-2) (\mu+\nu)^{1-m}\mu|\p_\mu W|^2\right) d\mu d\nu d\vol_{h_0} + \\
\ha \int_{\Omega_{a_0b}\times \p X} (\mu+\nu)^{-m}(\phi |\nabla_{h_0} W|^2 + B|W|^2) d\mu d\nu d\vol_{h_0} -  \\
\int_{\Omega_{a_0b} \times \p X} \mu (\nu-\mu)(\mu+\nu)^{-1-m}( |\phi'| |\nabla_{h_0} W|^2 + |B'||W|^2) d\mu d\nu d\vol_{h_0} \leq \\ T( ||q_1||_{1,\frac{m}{2}}+ ||q_2||_{1,\frac{m}{2}}) 
+ \int_{\Omega_{a_0b} \times \p X} G\mu (\mu+\nu)^{-m} (\p_\mu-\p_\nu) W d\mu d\nu d\vol_{h_0}.
\end{gathered} \label{intm}
\end{gather}
When $m=2$ this gives
\begin{gather*}
\ha \int_{\Omega_{a_0b}\times \p X}\left[ |\p_\nu W|^2 +  (\mu+\nu)^{-2}(\phi |\nabla_{h_0} W|^2 + B|W|^2) 
\right] d\mu d\nu d\vol_{h_0} -  \\
\int_{\Omega_{a_0b} \times \p X} (\mu+\nu)^{-1}( |\phi'| |\nabla_{h_0} W|^2 + |B'||W|^2) d\mu d\nu d\vol_{h_0} \leq \\ T( ||q_1||_{1,1}+ ||q_2||_{1,1}) 
+ T\int_{\Omega_{a_0b} \times \p X} |G|(\mu+\nu)^{-2} (\p_\mu-\p_\nu) W d\mu d\nu d\vol_{h_0}.
\end{gather*}

If $n>3,$ then $B>0$ and if $T$ is small, the second integral can be absorbed into the first.  When $n=3$ 
the argument used above shows that this can also be done. So we obtain
\begin{gather*}
\oq \int_{\Omega_{a_0b}\times \p X}\left[ |\p_\nu W|^2 +  (\mu+\nu)^{-2}(\phi |\nabla_{h_0} W|^2 + B|W|^2) 
\right] d\mu d\nu d\vol_{h_0} 
\leq \\ T( ||q_1||_{1,1}+ ||q_2||_{1,1}) 
+ T\int_{\Omega_{a_0b} \times \p X} |G|(\mu+\nu)^{-2} |(\p_\mu-\p_\nu) W| d\mu d\nu d\vol_{h_0}.
\end{gather*}

Now we repeat this argument by dropping the second integral in \eqref{divth} and integrate the remaining terms in $b$ with $a_0 \leq b \leq T.$

\begin{gather*}
\oq \int_{\Omega_{a_0T}\times \p X}\left[ |\p_\mu W|^2 +  (\mu+\nu)^{-2}(\phi |\nabla_{h_0} W|^2 + B|W|^2)  \right] d\mu d\nu d\vol_{h_0} 
\leq \\ T( ||q_1||_{1,1}+ ||q_2||_{1,1}) 
+ T\int_{\Omega_{a_0T} \times \p X} |G|(\mu+\nu)^{-2} |(\p_\mu-\p_\nu) W| d\mu d\nu d\vol_{h_0}.
\end{gather*}
We now add these estimates 
\begin{gather*}
\oq \int_{\Omega_{a_0T}\times \p X}\left[ |\p_\nu W|^2 +  |\p_\mu W|^2+ (\mu+\nu)^{-2}(\phi |\nabla_{h_0} W|^2 + B|W|^2) 
\right] d\mu d\nu d\vol_{h_0} 
\leq \\ 2 T( ||q_1||_{1,1}+ ||q_2||_{1,1}) 
+ 2 T\int_{\Omega_{a_0T} \times \p X} |G|(\mu+\nu)^{-2} (\p_\mu-\p_\nu) W d\mu d\nu d\vol_{h_0}.
\end{gather*}
By the Cauchy-Schwartz inequality 
\begin{gather*}
\int_{\Omega_{a_0T} \times \p X} |G|(\mu+\nu)^{-2} (\p_\mu-\p_\nu) W d\mu d\nu d\vol_{h_0}
\leq 2 |||G|||_{0,2} + 4 \int_{\Omega_{a_0T}}( |\p_\mu W|^2 +  |\p_\nu W|^2) d\mu d\nu d\vol_{h_0}
\end{gather*}
If $T$ is small this gives 
\begin{gather*}
 \int_{\Omega_{a_0T}\times \p X}\left[ |\p_\nu W|^2 +  |\p_\mu W|^2+ (\mu+\nu)^{-2}(\phi |\nabla_{h_0} W|^2 + B|W|^2) 
\right] d\mu d\nu d\vol_{h_0} 
\leq \\  C T( ||q_1||_{1,1}+ ||q_2||_{1,1} +|||G|||_{0,2}).
\end{gather*}
Equation \eqref{mainest2} then follows from \eqref{divth}.

 When $m=3$ equation \eqref{intm} gives
 \begin{gather*}
 \ha \int_{\Omega_{a_0b}\times \p X}(\mu+\nu)^{-2}(\mu |\p_\mu W|^2 + \nu|\p_\nu W|^2) d\mu d\nu d\vol_{h_0} + \\
\ha \int_{\Omega_{a_0b}\times \p X} (\mu+\nu)^{-3}(\phi |\nabla_{h_0} W|^2 + B|W|^2) d\mu d\nu d\vol_{h_0} -  \\
\int_{\Omega_{a_0b} \times \p X} (\mu+\nu)^{-2}( |\phi'| |\nabla_{h_0} W|^2 + |B'||W|^2) d\mu d\nu d\vol_{h_0} \leq \\ T( ||q_1||_{1,\tha}+ ||q_2||_{1,\tha})
+ \int_{\Omega_{a_0b} \times \p X} |G\mu (\mu+\nu)^{-3} (\p_\mu-\p_\nu) W| d\mu d\nu d\vol_{h_0}.
\end{gather*}

The Cauchy-Schwartz inequality gives
\begin{gather*}
\int_{\Omega_{a_0b} \times \p X} |G\mu (\mu+\nu)^{-3} (\p_\mu-\p_\nu) W |d\mu d\nu d\vol_{h_0}\leq
16 T |||G|||_{0,2}^2 + \frac{1}{4} \int_{\Omega_{a_0b}} \mu (\mu+\nu)^{-2}(|\p_\mu W|^2 +|\p_\nu W|^2) d\mu d\nu d\vol_{h_0} \leq \\
16 T |||G|||_{0,2}^2 + \frac{1}{4} \int_{\Omega_{a_0b}}  (\mu+\nu)^{-2}(\mu |\p_\mu W|^2 +\nu |\p_\nu W|^2) d\mu d\nu d\vol_{h_0} \leq \\
\end{gather*}
Thus we obtain
 \begin{gather*}
\oq \int_{\Omega_{a_0b}\times \p X}(\mu+\nu)^{-2}(\mu |\p_\mu W|^2 + \nu|\p_\nu W|^2) d\mu d\nu d\vol_{h_0} + \\
\ha \int_{\Omega_{a_0b}\times \p X} (\mu+\nu)^{-3}(\phi |\nabla_{h_0} W|^2 + B|W|^2) d\mu d\nu d\vol_{h_0} -  \\
\int_{\Omega_{a_0b} \times \p X} (\mu+\nu)^{-2}( |\phi'| |\nabla_{h_0} W|^2 + |B'||W|^2) d\mu d\nu d\vol_{h_0} \leq \\ T( ||q_1||_{1,\tha}+ ||q_2||_{1,\tha} +||| G|||_{0,2}) 
\end{gather*}
 For small $T_0,$ the third integral can be absorbed into the second and we obtain \eqref{mainest3}.  

Now we consider the case when $m=4.$  We deduce from \eqref{intm} that
\begin{gather}
\begin{gathered}
\ha \int_{\Omega_{a_0b}\times \p X}\left[ (\nu+\mu)^{-3}( (\nu-\mu)|\p_\nu W|^2 + 2 \mu|\p_\mu W|^2\right) d\mu d\nu d\vol_{h_0} + \\
\ha \int_{\Omega_{a_0b}\times \p X} (\mu+\nu)^{-4}(\phi |\nabla_{h_0} W|^2 + B|W|^2) d\mu d\nu d\vol_{h_0} -  \\
\int_{\Omega_{a_0b} \times \p X} \mu (\nu-\mu)(\mu+\nu)^{-5}( |\phi'| |\nabla_{h_0} W|^2 + |B'||W|^2) d\mu d\nu d\vol_{h_0} \leq \\ T( ||q_1||_{1,2}+ ||q_2||_{1,2}) 
+ \int_{\Omega_{a_0b} \times \p X} |G|\mu (\mu+\nu)^{-4} |(\p_\mu-\p_\nu) W| d\mu d\nu d\vol_{h_0}.
\end{gathered} \label{int4}
\end{gather}

We deduce from the Cauchy-Scwhartz inequality that
\begin{gather*}
 \int_{\Omega_{a_0b} \times \p X} |G|\mu (\mu+\nu)^{-4} ||(\p_\mu-\p_\nu) W d\mu d\nu d\vol_{h_0}
 \leq	\\ 64 |||G|||_{0,\fha}^2 + \frac{1}{4} \int_{\Omega_{a_0b}} \mu^2(\mu+\nu)^{-3} (|\p_\mu W|^2+|\p_\nu W|^2) d\mu d\nu d\vol_{h_0} \leq \\
 64 |||G|||_{0,\fha}^2 + \frac{1}{4} \int_{\Omega_{a_0b}} (\mu+\nu)^{-2} (\mu |\p_\mu W|^2+\nu |\p_\nu W|^2) d\mu d\nu d\vol_{h_0}, 
 \end{gather*}
and from \eqref{mainest3} we get that
 \begin{gather*}
  \int_{\Omega_{a_0b} \times \p X} |G|\mu (\mu+\nu)^{-4} ||(\p_\mu-\p_\nu) W| d\mu d\nu d\vol_{h_0}\leq
 C |||G|||_{0,\fha} + CT( ||q_1||_{1,2}+ ||q_2||_{1,2} +||| G|||_{0,2}). 
  \end{gather*}
 By substituting this into \eqref{int4} we find that
\begin{gather}
\begin{gathered}
\ha \int_{\Omega_{a_0b}\times \p X}\left[ (\nu+\mu)^{-3}( (\nu-\mu)|\p_\nu W|^2 + 2 \mu|\p_\mu W|^2\right) d\mu d\nu d\vol_{h_0} + \\
\ha \int_{\Omega_{a_0b}\times \p X} (\mu+\nu)^{-4}(\phi |\nabla_{h_0} W|^2 + B|W|^2) d\mu d\nu d\vol_{h_0} -  \\
\int_{\Omega_{a_0b} \times \p X} \mu (\nu-\mu)(\mu+\nu)^{-5}( |\phi'| |\nabla_{h_0} W|^2 + |B'||W|^2) d\mu d\nu d\vol_{h_0} \leq \\ T( ||q_1||_{1,2}+ ||q_2||_{1,2} + |||G|||_{2,0}) + |||G|||_{0,\fha}
\end{gathered} \label{int41}
\end{gather} 

 Again, when $T$ is small the third integral can be absorbed into the second and in particular 
 \eqref{mainest} follows from \eqref{int41}.  This proves Lemma \ref{trueestimates}.
 \end{proof}

Now we are ready to prove  Theorem \ref{smooth}.  
\begin{proof}
We will concentrate on the cases not covered by Proposition \ref{deg1}. So we assume that  if $n=3,$ 
$\int_{\p X} w \; d\vol_{h_0}=0.$

We can then apply Lemma \ref{trueestimates} to equation
\eqref{nonzero} with $G=0.$  We get that the solution $w$ to equation \eqref{weq3} satisfies

\begin{gather}
|||w|||_{0,2} \leq C(\widetilde{f_1},\widetilde{f_2}). \label{prim1}
\end{gather}
Since $\Delta_{h_0}$ commutes with the equation, we also have
\begin{gather}
|||\Delta_{h_0}^k w|||_{0,2} \leq C(\widetilde{f_1},\widetilde{f_2}), \;\ k=0,1,.... \label{prim2}
\end{gather}
Now we differentiate equation \eqref{weq3} with respect to the vector field $\p_\mu-\p_\nu.$
Let $W_j=(\p_\mu-\p_\nu)^j w.$ We find that $W_1$ satisfies
\begin{gather*}
\left((\mu+\nu)^2\p_\mu\p_\nu- \phi(\frac{2\mu\nu}{\mu+\nu}) \Delta_{h_0} - B(\frac{2\mu\nu}{\mu+\nu})\right)W_1=  2 (\nu-\mu)(\mu+\nu)^{-1}( \phi' \Delta_{h_0} w-B'w) \text{ in } \Omega_T\times \p X \\
W_1(\mu,\mu)=Q_1(\mu), \;\ \p_\mu W_1(\mu,\mu)=Z_1(\mu), \;\  
\end{gather*}
where $Q_1, Z_1 \in C^\infty([0,T]),$ depend on $\widetilde{f_1}$ and $\widetilde{f_2}$ and satisfy
\begin{gather*}
\p_\mu^k Q_1(0)=0, \;\  \p_\mu^k Z_1(0)=0, \;\  \;\ k=0,1,2,... 
\end{gather*}
Let $G_1= 2 (\nu-\mu)(\mu+\nu)^{-1}( \phi' \Delta_{h_0} w-B'w).$ In view of \eqref{prim1} and \eqref{prim2},
$|||G_1|||_{0,2} \leq C(\widetilde{f_1},\widetilde{f_2}).$ Then Lemma \ref{trueestimates} guarantees that
$|||W_1|||_{0,\tha} \leq C(\widetilde{f_1},\widetilde{f_2}).$  So Lemma \ref{decaym} implies that

\begin{gather}
|||w|||_{0,\fha} \leq C(\widetilde{f_1},\widetilde{f_2}), \;\ |||\Delta_{h_0}^k w|||_{0,\fha} \leq C(\widetilde{f_1},\widetilde{f_2}), \;\\ k=0,1,2... \label{prim21}
\end{gather}
But then $|||G_1|||_{0,\fha} \leq C(\widetilde{f_1},\widetilde{f_2})$ and Lemma \ref{trueestimates} guarantees that 
\begin{gather}
|||W_1|||_{0,2} \leq C(\widetilde{f_1},\widetilde{f_2}).\label{regw1}
\end{gather}  
Then Lemma \ref{decaym} gives that
\begin{gather}
|||w|||_{0,3} \leq C(\widetilde{f_1},\widetilde{f_2}), \;\ |||\Delta_{h_0} w|||_{0,3} \leq C(\widetilde{f_1},\widetilde{f_2}). \label{prim210}
\end{gather}
Now we differentiate \eqref{weq3} again with respect to $\p_\mu-\p_\nu.$ We find that
\begin{gather*}
\left((\mu+\nu)^2\p_\mu\p_\nu- \phi(\frac{2\mu\nu}{\mu+\nu}) \Delta_{h_0} - B(\frac{2\mu\nu}{\mu+\nu})\right)W_2=  2 (\nu-\mu)(\mu+\nu)^{-1}( \phi' \Delta_{h_0} -B') W_1 + \\
\left[-4(\mu+\nu)^{-1}( \phi' \Delta_{h_0} -B') + 4(\mu-\nu)^2(\mu+\nu)^{-2}(\phi'' \Delta_{h_0}-B'')\right]w
 \text{ in } \Omega_T\times \p X \\
W_2(\mu,\mu)=Q_2(\mu), \;\ \p_\mu W_1(\mu,\mu)=Z_2(\mu), \;\  
\end{gather*}
Let 
\begin{gather*}
G_2=2 (\nu-\mu)(\mu+\nu)^{-1}( \phi' \Delta_{h_0} -B') W_1 + \\
\left[-4(\mu+\nu)^{-1}( \phi' \Delta_{h_0} -B') + 4(\mu-\nu)^2(\mu+\nu)^{-2}(\phi'' \Delta_{h_0}-B'')\right]w
\end{gather*}
\end{proof}

In view of \eqref{regw1} and \eqref{prim210}, $|||G_2|||_{0,2} \leq C(\widetilde{f_1},\widetilde{f_2}),$ but then it follows from
Lemma \ref{trueestimates} that $|||W_2|||_{0,\tha} \leq C(\widetilde{f_1},\widetilde{f_2}).$ Then
Lemma \ref{decaym} implies that  $|||W_1|||_{0,\fha} \leq C(\widetilde{f_1},\widetilde{f_2}),$
and so  $|||w|||_{0,\frac{7}{2}} \leq C(\widetilde{f_1},\widetilde{f_2}).$ This implies that in fact
$|||G_2|||_{0,\fha} \leq C(\widetilde{f_1},\widetilde{f_2}),$ and therefore
$ |||W_2|||_{0,2}\leq C(\widetilde{f_1},\widetilde{f_2}).$ Now we differentiate the equation again
and repeat the argument.   We find that
\begin{gather*}
 ||| (\p_\mu-\p_\nu)^j w|||_{0,\fha} \leq  C(\widetilde{f_1},\widetilde{f_2}), \;\ j=1,2,... 
 \end{gather*}
and by Lemma \ref{decaym} we conclude that
\begin{gather}
 ||| w|||_{0,j} \leq  C(\widetilde{f_1},\widetilde{f_2}), \;\ 
  ||| \Delta_{h_0}^k w|||_{0,j} \leq  C(\widetilde{f_1},\widetilde{f_2}), \;\ k, j=1,2,...  \label{final1}
 \end{gather}
 Now we go back to equation \eqref{divth} and apply it to the solution $w$ to
 \eqref{weq3}.  In this case $G=0$ and we obtain
\begin{gather}
\begin{gathered}
\ha \int_{\Sigma_1\times \p X}\left[ (\mu+\nu)^{2-m} |\p_\nu w|^2 +(\mu+\nu)^{-m}(\phi |\nabla_{h_0} w|^2 + B|w|^2)\right] d\nu d\vol_{h_0} + \\
\ha \int_{\Sigma_2 \times \p X}\left[ (\mu+\nu)^{2-m} |\p_\mu w|^2 +
(\mu+\nu)^{-m}(\phi |\nabla_{h_0} w|^2 + B|W|^2)\right] d\mu  d\vol_{h_0}  + \\
\frac{m-2}{2}\int_{\Omega_{ab} \times \p X} (\mu+\nu)^{1-m}((\p_\mu w)^2 -(\p_\nu w)^2) d\mu d\nu 
d\vol_{h_0} + \\
\int_{\Omega_{ab} \times \p X} (\nu-\mu)(\mu+\nu)^{-1-m}( \phi'|\nabla_{h_0} w|^2 + B'|w|^2) d\mu d\nu d\vol_{h_0}= \\
\frac{1}{2\sqrt{2}} \int_{\Sigma_3 \times \p X} (2\mu)^{2-m} ((\widetilde{f_2})^2  +
(\widetilde{f_2}- \widetilde{f_1})^2) + (2\mu)^{-m} (\phi(\mu)|\nabla_{h_0} \widetilde{f_1}
|^2+ B(\mu) |\widetilde{f_1}|^2) \; d\mu
d\vol_{h_0}.
\end{gathered}\label{divthfinal}
\end{gather} 
 
    The term
 \begin{gather*}
\int_{\Omega_{ab} \times \p X} (\mu+\nu)^{1-m}((\p_\mu w)^2 -(\p_\nu w)^2) d\mu d\nu 
d\vol_{h_0} = \int_{\Omega_{ab} \times \p X} (\mu+\nu)^{1-m}
((\p_\mu  -\p_\nu) w)((\p_\mu  +\p_\nu) w)d\mu d\nu d\vol_{h_0}\leq \\
2 \int_{\Omega_{ab} \times \p X} (\mu+\nu)^{2-2m}
|(\p_\mu  -\p_\nu) w|^2d\mu d\nu d\vol_{h_0}+
2 \int_{\Omega_{ab} \times \p X} 
|(\p_\mu  +\p_\nu) w|^2d\mu d\nu d\vol_{h_0}\leq C(\widetilde{f_1},\widetilde{f_2}), 
\end{gather*}
and from \eqref{final1} we have that
\begin{gather*}
\left|\int_{\Omega_{ab} \times \p X} (\nu-\mu)(\mu+\nu)^{-1-m}( \phi'|\nabla_{h_0} w|^2 + B'|w|^2) d\mu d\nu d\vol_{h_0}\right| \leq  C(\widetilde{f_1},\widetilde{f_2}).
\end{gather*}
Therefore we conclude that
\begin{gather}
\begin{gathered}
\ha \int_{\Sigma_1\times \p X}\left[ (\mu+\nu)^{2-m} |\p_\nu w|^2 +(\mu+\nu)^{-m}(\phi |\nabla_{h_0} w|^2 + B|w|^2)\right] d\nu d\vol_{h_0} + \\
\ha \int_{\Sigma_2 \times \p X}\left[ (\mu+\nu)^{2-m} |\p_\mu w|^2 +
(\mu+\nu)^{-m}(\phi |\nabla_{h_0} w|^2 + B|W|^2)\right] d\mu  d\vol_{h_0}  \leq
C(\widetilde{f_1},\widetilde{f_2}).
\end{gathered}\label{divthfinal1}
\end{gather} 
Now we write, for $m>2,$
\begin{gather}
\begin{gathered}
|(\mu+\nu)^{1-m} w(\mu,\nu)| \leq 2\mu^{1-m} |w(\mu,\mu)|+ 2(\mu+\nu)^{2-m}\int_{\mu}^\nu
|\p_s w(\mu,s)|^2 ds \leq \\
2\mu^{1-m} |w(\mu,\mu)|+ 2\int_{\mu}^\nu (\mu+s)^{2-m} |\p_s w(\mu,s)|^2 ds
\leq C(\widetilde{f_1},\widetilde{f_2}).
\end{gathered}\label{important}
\end{gather}
Now we argue as we did in the case of $w_0,$ and conclude that 
$w\in C^\infty(\overline{\Omega_T} \times \p X).$
This ends the proof of Theorem \ref{smooth}.
\end{proof}

To prove Theorem \ref{main} we need the following:
\begin{lemma}\label{vanishing} Suppose that $f \in \mcs(X)$ and  $\mcr_+(0,f)(s,y)=0$
for $s<s_0<<0.$  
Let $u$ be the solution to \eqref{weqnew} and let $w$ be the function defined by $u$ in
\eqref{wfu}.  Then 
\begin{gather}
\begin{gathered}
\p_{\mu}^k w(0,\nu,y)=0 \text{ if } \nu\leq -\frac{1}{s_0}, \;\ k=0,1,...  \\
\p_{\nu}^k w(\mu,0,y)=0 \text{ if } \mu\leq -\frac{1}{s_0},  \;\ k=0,1,... 
\end{gathered}\label{van1}
\end{gather}
\end{lemma}
\begin{proof}
Suppose $\mcr_+(0,f)(s,y)=0$ for $s<s_0.$  Since the initial data is of the form $(0,f),$ the solution $u$ of
\eqref{weqnew} is odd in time, and therefore $\mcr_-(0,f)(s,y)=0$ if $s>-s_0.$  Since
$w$ is smooth up to $\{\mu=0\}$ and $\{\nu=0\},$  this means that
\begin{gather*}
 \p_\nu w(0,\nu,y)=0 \;\ \text{ if } \nu\leq -\frac{1}{s_0}, \;\ 
 \p_\mu w(\mu,0,y)=0 \;\ \text{ if } \mu\leq -\frac{1}{s_0}.
\end{gather*}
We know from \eqref{important} that $w$ vanishes to infinite order at $\{\mu=\nu=0\}.$ Therefore
\begin{gather*}
w(0,\nu,y)=0 \;\ \text{ if } \nu\leq -\frac{1}{s_0}, \;\
 w(\mu,0,y)=0 \;\ \text{ if } \mu\leq -\frac{1}{s_0}.
\end{gather*} 

From the equation \eqref{weq3},  we conclude that 
\begin{gather*}
\text{if } \mu\not=0, \;\ \p_\mu \p_\nu w(\mu,0,y)= \mu^{-2}( \phi(0) \Delta_{h_0} + B(0)) w(\mu,0)=0 \\
\text{if } \nu\not=0, \;\ \p_\mu \p_\nu w(0,\nu,y)= \nu^{-2}( \phi(0) \Delta_{h_0} + B(0)) w(0,\nu)=0.
\end{gather*}
Using that $w$ vanishes to infinite order at $\{\mu=\nu=0\},$
\begin{gather*}
\p_\mu w(0,\nu,y)=0 \;\ \text{ if } \nu\leq -\frac{1}{s_0}, \\
\p_\nu w(\mu,0,y)=0 \;\ \text{ if } \mu\leq -\frac{1}{s_0}.
\end{gather*} 
 
    Now we differentiate equation \eqref{weq3} with respect to
 $\mu$ we find that
 \begin{gather*}
 (\mu+\nu)^2 \p_\mu^2\p_\nu w+ \phi \Delta_{h_0} \p_\mu w+ B \p_\mu w+ 2(\mu+\nu) \p_\mu\p_\nu w +
 \nu^2(\mu+\nu)^{-2} \Delta_{h_0} w+ \nu^2(\mu+\nu)^{-2} B' w=0
 \end{gather*}
 Hence, 
 \begin{gather*}
\text{ if } \nu\not=0,  \;\   \p_\mu^2\p_\nu w(0,\nu,y)= 0.
\end{gather*}
Since $w$ is smooth, and vanishes to infinite order at $\mu=\nu=0,$
\begin{gather*}
\p_\mu^2 w(0,\nu,y)=0 \;\ \text{ if } \nu\leq -\frac{1}{s_0}.
\end{gather*} 
By symmetry,
\begin{gather*}
\p_\nu^2 w(\mu,0,y)=0 \;\ \text{ if } \mu\leq -\frac{1}{s_0}.
\end{gather*} 
Since away from $\{\mu=\nu=0\}$ the coefficients of \eqref{weq3} are smooth,
we can repeat the argument to prove that all derivatives of $w$ vanish at
$\{\mu=0, \;\ \nu\leq -\frac{1}{s_0} \}\cup \{\nu=0, \;\ \mu\leq -\frac{1}{s_0}\}.$
\end{proof}
\section{Carleman Estimates}

Let $w$ be a solution to  $\eqref{weq3}$  with $\widetilde{f_1}=0$ and 
$\widetilde{f_2}=\widetilde{f}=\frac{1}{2\mu^2} F(\mu)f(\mu,y),$ and let 
\begin{gather*}
w_k(\mu,\nu)=\lan w(\mu,\nu,y), \phi_k(y) \ran_{L^2(\p X,d\vol_{h_0})},
\end{gather*}
where $\phi_k,$ $k=1,2,...$ are the eigenfunctions of $\Delta_{h_0}$ with eigenvalue
$\lambda_k,$ where $0=\la_1<\lambda_2,$ and $\la_k\leq \lambda_{k+1},$ for $k>1,$ and
$\la_k \rightarrow \infty.$ Then
 $w_k\in C^\infty([0,T]\times [0,T])$ satisfies
\begin{gather*}
\left((\mu+\nu)^2\p_\mu\p_\nu + F_k(\mu,\nu)\right) w_k=0 \text{ in } [0,T]\times [0,T], \\
w_k(\mu,\mu)=0, \;\ \p_\mu w_k(\mu,\mu)=\widetilde{f_{k}}(\mu)= \lan \widetilde{f}, \phi_k\ran, \;\ j=1,2.
\end{gather*}
Here $F_k(\mu,\nu)=\la_k \phi(\mu,\nu) + B(\mu,\nu).$

By assumption $w$ vanishes to infinite order at $\{\mu=0\}$ and $\{\nu=0\}.$ Thus so does
$w_k(\mu,\nu),$  $k=1,2,...$  We will prove that under these assumptions there exists $\eps>0,$ independent of $k$
 such that $w_k(\mu,\nu)=0$ if $\mu\leq \eps$ and $\nu\leq \eps.$  In particular
 $\widetilde{f_k}(\mu) =0,$  $k=0,1,2,...$ if $\mu\leq \eps,$ and hence $f(x)=0$ if $x\leq \eps.$
 
 It is convenient to work with coordinates $r$ and $\tau$ defined in \eqref{rotate}.
 The main ingredient of the proof is
 \begin{lemma}\label{carleman} Let  $U$ be a neighborhood of $(0,0).$ Let
 \begin{gather*}
P_{k}=r^2\p_r^2-r^2\p_\tau^2 + F_k(\tau,r),
\end{gather*}
where $F_k(r,\tau)= \la_k \phi(\frac{r^2-\tau^2}{2r}) -B(\frac{r^2-\tau^2}{2r}).$
Then there exists  $C>0,$  independent of $k,$ and $\gamma_0=\gamma_0(k)$ such that for every 
$\gamma>\gamma_0,$ and every $u \in C_0^\infty(U),$ which is supported
in $\{(\tau, r): \;\ r\geq 0 \text{ and } -r \leq \tau \leq r\},$
\begin{gather}
||r^{-\gamma-2} P_k u||^2 \geq  C\left( \gamma^2 ||r^{-1}\p_r r^{-\gamma} u||^2 +
\gamma^2 ||r^{-1-\gamma} \p_{\tau} u||^2 +
\gamma^4 || r^{-\gamma-2} u||^2\right). \label{carle1}
\end{gather}
Here $||\cdot ||$ denotes the $L^2(U)$ norm.
\end{lemma}
\begin{proof} 
Let 
\begin{gather*}
P_{\gamma,k}= r^{-\gamma-2} P_k r^\gamma, \;\ P_{\gamma,k}= r^{-2} P_k + 2\gamma r^{-1} \p_r + \gamma(\gamma-1) r^{-2}.
\end{gather*}

The support of $u$ is contained in $\{r\geq 0\}$ and $\{ r\geq \tau \geq -r\}.$ So we write
\begin{gather*}
u=r^{\gamma} v, \text{  and } r^{-\gamma-2}P_k u=P_{\gamma,k} v.
\end{gather*}
 We have
\begin{gather}
\begin{gathered}
||P_{\gamma,k} v||^2= ||r^{-2} P_k v||^2 + 4\gamma^2|| r^{-1}\p_r v||^2 +\gamma^2(\gamma-1)^2 ||r^{-2} v||^2
+ 4\gamma\lan r^{-2} P_k v, r^{-1} \p_r v\ran + \\
2\gamma(\gamma-1) \lan r^{-2} P_k v, r^{-2} v\ran  + 4\gamma^2(\gamma-1)\lan r^{-1}\p_r v, r^{-2} v\ran.
\end{gathered} \label{carle2}
\end{gather}
Now we integrate by parts to compute the inner products
$\lan r^{-2}  P_k v, r^{1} \p_r \ran,$
$\lan r^{-2} P_k v, r^{-2} v\ran$ and $\lan r^{-1}\p_r v, r^{-2} v\ran.$
We begin with 
\begin{gather*}
\lan r^{-1}\p_r v, r^{-2} v\ran=  \ha \int r^{-3} \p_r v^2 dr d\tau=\tha ||r^{-2} v||^2. 
\end{gather*}
$\lan r^{-2} P_k v, r^{-1} \p_r v\ran$ has two terms:
\begin{gather*}
\lan (\p_r^2 -\p_t^2) v, r^{-1} \p_r v\ran= \ha \int r^{-1} \p_r(\p_r v)^2 dr d\tau -
\int r^{-1} \p_t^2 \p_r v dr d\tau= \ha ||r^{-1} \p_r v||^2 + \ha ||r^{-1} \p_t v||^2,
\end{gather*}
and
\begin{gather*}
\lan   r^{-2} F_k v, r^{-1} \p_r v\ran= \ha \int r^{-3} F_k \p_r v^2 dr d\tau=
\tha \int r^{-4}( F_k +r\p_r F_k) v^2 dr d\tau \geq - C(\la_k+1)||r^{-2} v||^2.
\end{gather*}
$\lan r^{-2} P_k v, r^{-1} \p_r v\ran$ also has two terms:
\begin{gather*}
\lan (\p_r^2 -\p_t^2) v , r^{-2} \ran= 3||r^{-2} v||^2 -||r^{-1} \p_r v||^2 +||r^{-1} \p_t v||^2.
\end{gather*}
and
\begin{gather*}
\lan r^{-2} F_k v, r^{-2} v\ran \geq -C(\la_k +1) ||r^{-2}v||^2.
\end{gather*}

Putting these estimates together we find that
\begin{gather*}
||P_{\gamma,k} v||^2 \geq (2\ga^2+4\ga) ||r^{-1} \p_r v||^2 + 2\ga^2 ||r^{-1} \p_t v||^2 + \\
(\ga^2(\ga-1)^2 + 6\ga^2(\ga-1) + 6\ga(\ga-1) -4C(\la_k+1)\ga-2\ga(\ga-1)\la_k) ||r^{-2}v||^2.
\end{gather*}
Thus, if $\ga_0>> C(\la_k+1),$
\begin{gather*}
||P_{\gamma,k} v||^2 \geq C(\ga^2 ||r^{-1} \p_r v||^2 + \ga^2 ||r^{-1} \p_t v||^2 + 
\ga^4 ||r^{-2}v||^2).
\end{gather*}

Since $v=r^{-\ga} u$ this implies that 
\begin{gather}
||r^{-\ga-2} P_k u||^2 \geq C(\ga^2  ||r^{-\ga} \p_r r^{-1} u||^2+ \ga^2||r^{-\ga-1} \p_t u||^2 +
\ga^4 ||r^{-\ga-2} u||^2, \;\ \ga > \ga_0. \label{carle3}
\end{gather}
\section{ Proof of Theorem \ref{main}}
To prove Theorem \ref{main} we apply Lemma \ref{carleman} to $\chi(r) w_k,$ with
$\chi \in C^\infty(-T,T),$  such that $\chi(r)=1$ if $|r|<\frac{T}{4}$ and $\chi(r)=0$ if $|r|> \frac{T}{2}.$

As a consequence of  \eqref{carle3},
\begin{gather*}
|| r^{-\ga-2} P_k \chi w_k|| \geq C \ga^4 ||r^{-\ga-2} \chi w_k||^2.
\end{gather*}

Notice that, since $P_k w_k=0,$
\begin{gather*}
P_k \chi(r) w_k=r^2 \chi''(r) w_k + 2r^2\chi'(r) \p_r w_k,
\end{gather*}
and therefore $P_k \chi(r) w_k$ is supported in $\frac{T}{4} \leq r \leq \frac{T}{2}.$  Thus
\begin{gather*}
||r^{-\ga-2} P_k \chi w_k ||^2 \leq  C(\chi,w_k) T^{-2\ga-2} 4^{2\ga+2}
\end{gather*}
But since $\chi(r)=1$ if $r<T/4,$
\begin{gather*}
||r^{-\ga-2} \chi w_k || \geq ||r^{-\ga-2}  w_k ||_{L^2(B(0,\frac{T}{4}))} \geq T^{-2\ga-2} 4^{2\ga+2}
 || w_k ||_{L^2(B(0,\frac{T}{4}))},
\end{gather*}
where $ || w_k ||_{L^2(B(0,\frac{T}{4}))}$ is the $L^2$ norm of
$ w_k $ on the ball centered at $(0,0)$ and radius  $\frac{T}{4}.$
So we have
\begin{gather*}
C(\chi,w_k) \geq C \gamma^4 || w_k ||_{L^2(B(0,\frac{T}{4}))}, \;\ \gamma \geq \gamma_0.
\end{gather*}
Letting $\gamma\rightarrow \infty$ gives
\begin{gather*}
 ||w_k ||_{L^2(B(0,\frac{T}{4}))}=0.
 \end{gather*}
 Since $T$ does not depend on $k$ we conclude that $w=0$ in $B(0,\frac{T}{4}).$  In particular this shows that $f(x)=0$ if $x\leq \frac{T}{4}.$ Therefore $f$ is compactly supported. Then Theorem \ref{main} follows from Theorem \ref{compsup}.

 Notice that the weight in the Carleman estimate depends on the eigenvalue $\lambda_k.$   The bigger the eigenvalue, the larger the parameter $\gamma$ has to be.

 \qed

\section{acknowledgements}

This research was partially funded by the NSF grant DMS 0500788.

This paper was written in the summer 2007 when I visited
the Mathematics  Department of  the Federal University of Pernambuco, Recife, Brazil.  I would like to thank Prof. Fernando Cardoso for the hospitality.  I would also like to thank the CNPq (Brazil) for the financial support during my stay in Recife.

\end{proof}

\end{document}